\documentclass{article}
\usepackage{amsmath}
\usepackage{amssymb}
\usepackage{amsfonts}

\newtheorem{theorem}{Theorem}

\newtheorem{lemma}[theorem]{Lemma}

\newenvironment{proof}[1][Proof]{\noindent \textbf{#1.} }{\  \rule{0.5em}{0.5em}}
\input{tcilatex}
\begin{document}

\author{Zafer \c{S}iar$^{a}$ and Refik Keskin$^{b}$ \and $\ ^{a}$Bing\"{o}l
University, Mathematics Department, Bing\"{o}l/TURKEY$\ \ $ \and $^{b}$%
Sakarya University, Mathematics Department, Sakarya/TURKEY$\ $ \and $^{a}$%
zsiar@bingol.edu.tr, $^{b}$rkeskin@sakarya.edu.tr,}
\title{Repdigits in $k-$generalized Pell sequence}
\maketitle

\begin{abstract}
Let $k\geq 2$ and let $(P_{n}^{(k)})_{n\geq 2-k}$ be $k$-generalized Pell
sequence defined by 
\begin{equation*}
P_{n}^{(k)}=2P_{n-1}^{(k)}+P_{n-2}^{(k)}+...+P_{n-k}^{(k)}
\end{equation*}%
for $n\geq 2$ with initial conditions 
\begin{equation*}
P_{-(k-2)}^{(k)}=P_{-(k-3)}^{(k)}=\cdot \cdot \cdot
=P_{-1}^{(k)}=P_{0}^{(k)}=0,P_{1}^{(k)}=1.
\end{equation*}%
In this paper, we deal with the Diophantine equation 
\begin{equation*}
P_{n}^{(k)}=d\left( \frac{10^{m}-1}{9}\right)
\end{equation*}%
in positive integers $n,m,d$ with $m\geq 2$ and $1\leq d\leq 9$. We show
that repdigits in the sequence $\left( P_{n}^{(k)}\right) _{n\geq 2-k}$ ,
which have at least two digits, are the numbers\ $P_{5}^{(3)}=33$ and $%
P_{6}^{(4)}=88.$
\end{abstract}



\section{\protect\bigskip Introduction}

Let $k\geq 2$ be an integer. Let the linear recurrence sequence $\left(
G_{n}^{(k)}\right) _{n\geq 2-k}$ of order $k$ define by 
\begin{equation}
G_{n}^{(k)}=rG_{n-1}^{(k)}+G_{n-2}^{(k)}+\ldots +G_{n-k}^{(k)}\text{ }
\label{0}
\end{equation}%
for $n\geq 2$ with the initial conditions $%
G_{-(k-2)}^{(k)}=G_{-(k-3)}^{(k)}=\cdots =G_{-1}^{(k)}=0,$ $G_{0}^{(k)}=a,$
and $G_{1}^{(k)}=b.$ For $(a,b,r)=(0,1,1)$ and $(a,b,r)=(2,1,1),~$the
sequence $\left( G_{n}^{(k)}\right) _{n\geq 2-k}$ is called $k-$generalized
Fibonacci sequence $\left( F_{n}^{(k)}\right) _{n\geq 2-k}$ and $k-$%
generalized Lucas sequence $\left( L_{n}^{(k)}\right) _{n\geq 2-k}$(see \cite%
{luca2, Luca1}). For $(a,b,r)=(0,1,2)$ and $(a,b,r)=(2,2,2),~$the sequence $%
\left( G_{n}^{(k)}\right) _{n\geq 2-k}$ is called $k-$generalized Pell
sequence $\left( P_{n}^{(k)}\right) _{n\geq 2-k}$ and $k-$generalized
Pell-Lucas sequence $\left( Q_{n}^{(k)}\right) _{n\geq 2-k},$ respectively
(see \cite{klc1}). The terms of these sequences are called $k-$generalized
Fibonacci numbers, $k-$generalized Lucas numbers, $k-$generalized Pell
numbers and $k-$generalized Pell-Lucas numbers, respectively. When $k=2,$ we
have Fibonacci, Lucas, Pell and Pell-Lucas sequences, $\left( F_{n}\right)
_{n\geq 0},~\left( L_{n}\right) _{n\geq 0},$ $\left( P_{n}\right) _{n\geq
0}, $ and $\left( Q_{n}\right) _{n\geq 0}$, respectively.

A repdigit is a positive integer whose digits are all equal. Recently, some
mathematicians have investigated the repdigits in the above sequences for $%
k=2$ or general $k.$ In \cite{Luca4}, Luca determined that the largest
repdigits in the sequences $\left( F_{n}^{(2)}\right) _{n\geq 0}$ and $%
\left( L_{n}^{(2)}\right) _{n\geq 0}$ are $F_{10}^{(2)}=55$ and $%
L_{5}^{(2)}=11.$ In \cite{Faye}, the authors have found all repdigits in the
sequences $\left( P_{n}^{(2)}\right) _{n\geq 0}$ and $\left(
Q_{n}^{(2)}\right) _{n\geq 0}$. Here, they showed that the largest repdigits
in these sequences are $P_{3}^{(2)}=5$ and $Q_{2}^{(2)}=6.$ In \cite{marqu},
Marques proved that the largest repdigits in the sequence $\left(
F_{n}^{(3)}\right) _{n\geq -1}$ are $F_{8}^{(3)}=44.$ Besides, for general
case of $k,$ in \cite{bravo}, Bravo and Luca handled the Diophantine
equation 
\begin{equation}
F_{n}^{(k)}=d\left( \frac{10^{m}-1}{9}\right)  \label{1}
\end{equation}%
and showed that this equation has the solutions $%
(n,k,d,m)=(10,2,5,2),(8,3,4,2)$ in positive integers $n,m,k,d$ with $k\geq
2,~1\leq d\leq 9$ and $m\geq 2$. The same authors, in \cite{Luca1},
considered the equation (\ref{1}) for the sequence $\left(
L_{n}^{(k)}\right) _{n\geq 2-k},$ and they have given the solutions of this
equation by $(n,k,d,m)=(5,2,1,2),(5,4,2,2).$

In this paper, we will deal with the Diophantine equation 
\begin{equation}
P_{n}^{(k)}=d\left( \frac{10^{m}-1}{9}\right)  \label{1.1}
\end{equation}%
in positive integers $n,m,d$ with $m\geq 2$ and $1\leq d\leq 9$. We will
show that the repdigits in the sequence $\left( P_{n}^{(k)}\right) _{n\geq
2-k}$ , which have at least two digits, are the numbers\ $P_{5}^{(3)}=33$
and $P_{6}^{(4)}=88.$

\section{Preliminaries}

It can be seen that the characteristic polynomial of the sequence $\left(
P_{n}^{(k)}\right) _{n\geq 2-k}$ is%
\begin{equation}
\Psi _{k}(x)=x^{k}-2x^{k-1}-\cdots -x-1.  \label{0.2}
\end{equation}%
We know from Lemma 1 in given \cite{zang} that this polynomial has exactly
one positive real root located between $2$ and $3.$ We denote the roots of
the polynomial in (\ref{0.2}) by $\alpha _{1},\alpha _{2},\ldots ,\alpha
_{k}.$ Particuarly, let $\alpha =\alpha _{1}$ denote positive real root of $%
\Psi _{k}(x)$. The positive real root $\alpha =\alpha (k)$ is called
dominant root of $\Psi _{k}(x)$. The other roots are strictly inside the
unit circle. In \cite{luca3}, the Binet- like formula for $k-$ generalized
Pell numbers are given by 
\begin{equation}
P_{n}^{(k)}=\dsum\limits_{j=1}^{k}\frac{(\alpha _{j}-1)}{\alpha
_{j}^{2}-1+k(\alpha _{j}^{2}-3\alpha _{j}+1)}\alpha _{j}^{n}.  \label{x1}
\end{equation}%
It was also showed in \cite{luca3} the contribution of the roots inside the
unit circle to the formula (\ref{0.2}) is very small, namely that the
approximati\i on%
\begin{equation}
\left\vert P_{n}^{(k)}-g_{k}(\alpha )\alpha ^{n}\right\vert <\frac{1}{2}
\label{w}
\end{equation}%
holds for all $n\geq 2-k$, where 
\begin{equation}
g_{k}(z)=\frac{z-1}{(k+1)z^{2}-3kz+k-1}.  \label{w1}
\end{equation}%
The proof of the following inequality is given in \cite{zfr}. 
\begin{equation}
\left\vert \frac{(\alpha _{j}-1)}{\alpha _{j}^{2}-1+k(\alpha
_{j}^{2}-3\alpha _{j}+1)}\right\vert \leq 2  \label{1.4}
\end{equation}%
for $k\geq 2,$ where $\alpha _{j}$'s for $j=1,2,\ldots ,k$ are the roots of
the polynomial in (\ref{0.2}).

Throughout this paper, $\alpha $ denotes the positive real root of the
polynomial given in (\ref{0.2}). The following relation between $\alpha $
and $P_{n}^{(k)}$ given in \cite{luca3} is valid for all $n\geq 1.$ 
\begin{equation}
\alpha ^{n-2}\leq P_{n}^{(k)}\leq \alpha ^{n-1}.  \label{1.3}
\end{equation}%
Also, K\i l\i \c{c} \cite{klc1} proved that 
\begin{equation}
P_{n}^{(k)}=F_{2n-1}  \label{1.2}
\end{equation}%
for all $1\leq n\leq k+1.$

\begin{lemma}
\label{L1}\emph{(\cite{luca3}, Lemma 3.2)}Let $k,l\geq 2$ be integers. Then

\emph{(a)} If $k>l,$ then $\alpha (k)>\alpha (l),$ where $\alpha (k)$ and $%
\alpha (l)$ are the values of $\alpha $ relative to $k$ and $l,$
respectively.

\emph{(b) }$\varphi ^{2}(1-\varphi ^{-k})<\alpha <\varphi ^{2},$ where $%
\varphi $ is golden ratio.

\emph{(c)}$~g_{k}(\varphi ^{2})=\frac{1}{\varphi +2}.$

\emph{(d) }$0.276<g_{k}(\alpha )<0.5.$
\end{lemma}

For solving the equation (\ref{1.1}), we use linear forms in logarithms and
Baker's theory. For this, we will give some notions, theorem, and lemmas
related to linear forms in logarithms and Baker's theory.

Let $\eta $ be an algebraic number of degree $d$ with minimal polynomial 
\begin{equation*}
a_{0}x^{d}+a_{1}x^{d-1}+\cdots +a_{d}=a_{0}\dprod\limits_{i=1}^{d}\left(
X-\eta ^{(a)}\right) \in \mathbb{Z}[x],
\end{equation*}%
where the $a_{i}$'s are integers with $\gcd (a_{0},\ldots ,a_{n})=1$ and $%
a_{0}>0$ and $\eta ^{(a)}$'s are conjugates of $\eta .$ Then 
\begin{equation}
h(\eta )=\frac{1}{d}\left( \log a_{0}+\dsum\limits_{i=1}^{d}\log \left( \max
\left\{ |\eta ^{(a)}|,1\right\} \right) \right)  \label{2.1}
\end{equation}%
is called the logarithmic height of $\eta .$ In particularly, if $\eta =a/b$
is a rational number with $\gcd (a,b)=1$ and $b\geq 1,$ then $h(\eta )=\log
\left( \max \left\{ |a|,b\right\} \right) .$

We give some properties of the logarithmic height whose proofs can be found
in \cite{yann}:

\begin{equation}
h(\eta \pm \gamma )\leq h(\eta )+h(\gamma )+\log 2,  \label{2.2}
\end{equation}%
\begin{equation}
h(\eta \gamma ^{\pm 1})\leq h(\eta )+h(\gamma ),  \label{2.3}
\end{equation}%
\begin{equation}
h(\eta ^{m})=|m|h(\eta ).  \label{2.4}
\end{equation}%
In \cite{zfr}, using the above properties of the logarithmic height, the
authors have proved the inequality 
\begin{equation}
h(g_{k}(\alpha ))\leq 4\log k\text{ for }k\geq 3,  \label{2.5}
\end{equation}%
which will be used in the main theorem, where $g_{k}(\alpha )$ is as defined
in (\ref{w1}). Now we give a theorem deduced from Corollary 2.3 of Matveev 
\cite{Mtv} and provides a large upper bound for the subscript $n$ in the
equation (\ref{1.1}) (also see Theorem 9.4 in \cite{Bgud}).

\begin{theorem}
\label{T2} Assume that $\gamma _{1},\gamma _{2},\ldots ,\gamma _{t}$ are
positive real algebraic numbers in a real algebraic number field $\mathbb{K}$
of degree $D$, $b_{1},b_{2},\ldots ,b_{t}$ are rational integers, and 
\begin{equation*}
\Lambda :=\gamma _{1}^{b_{1}}\cdots \gamma _{t}^{b_{t}}-1
\end{equation*}%
is not zero. Then 
\begin{equation*}
|\Lambda |>\exp \left( -1.4\cdot 30^{t+3}\cdot t^{4.5}\cdot D^{2}(1+\log
D)(1+\log B)A_{1}A_{2}\cdots A_{t}\right) ,
\end{equation*}%
where 
\begin{equation*}
B\geq \max \left\{ |b_{1}|,\ldots ,|b_{t}|\right\} ,
\end{equation*}%
and $A_{i}\geq \max \left\{ Dh(\gamma _{i}),|\log \gamma _{i}|,0.16\right\} $
for all $i=1,\ldots ,t.$
\end{theorem}

The following lemma was proved by Dujella and Peth\H{o} \cite{duj} and is a
variation of a lemma of Baker and Davenport \cite{Baker}. This lemma will be
used to reduce the upper bound for the subscript $n$ in the equation (\ref%
{1.1}). For any real number $x,$ we let $||x||=\min \left\{ |x-n|:n\in 
\mathbb{Z}
\right\} $ be the distance from $x$ to the nearest integer.

\begin{lemma}
\label{L2}Let $M$ be a positive integer, let $p/q$ be a convergent of the
continued fraction of the irrational number $\gamma $ such that $q>6M,$ and
let $A,B,\mu $ be some real numbers with $A>0$ and $B>1.$ Let $\epsilon
:=||\mu q||-M||\gamma q||.$ If $\epsilon >0,$ then there exists no solution
to the inequality 
\begin{equation*}
0<|u\gamma -v+\mu |<AB^{-w},
\end{equation*}%
in positive integers $u,v,$ and $w$ with 
\begin{equation*}
u\leq M\text{ and }w\geq \frac{\log (Aq/\epsilon )}{\log B}.
\end{equation*}
\end{lemma}

The following lemma can be found in \cite{weger}.

\begin{lemma}
\label{L3} Let $a,x\in 
\mathbb{R}
.$ If $0<a<1$ and $\left\vert x\right\vert <a,$ then 
\begin{equation*}
\left\vert \log (1+x)\right\vert <\frac{-\log (1-a)}{a}\cdot \left\vert
x\right\vert
\end{equation*}%
and 
\begin{equation*}
\left\vert x\right\vert <\frac{a}{1-e^{-a}}\cdot \left\vert
e^{x}-1\right\vert .
\end{equation*}
\end{lemma}

\section{Main Theorem}

\begin{theorem}
\label{T4}The only solution of Diophantine equation\emph{\ (\ref{1.1})} in
positive integers $(n,m)$ with $1\leq d\leq 9$ are given by $%
(n,k,d,m)=(5,3,3,2),(6,4,8,2).$
\end{theorem}

\proof%
Assume that $P_{n}^{(k)}=d\left( \frac{10^{m}-1}{9}\right) $ with $n\geq 1$, 
$m,k\geq 2$ and $1\leq d\leq 9.$ If $1\leq n\leq k+1,$ then we have $d\left( 
\frac{10^{m}-1}{9}\right) =P_{n}^{(k)}=F_{2n-1}$ by (\ref{1.2}). In this
case we get $n=1,2,3$ by Theorem $1$ given in \cite{Luca4}. But, these
values of $n$ yields to $m=1,$ a contradiction. Then we suppose that $n\geq
k+2.$ If $k=2,$ then $n\geq 4$ and we have $P_{n}=d\left( \frac{10^{m}-1}{9}%
\right) ,$ which implies that $n=0,1,2,3$ by Theorem 1.1 given in \cite{Faye}%
. Again, since $m\geq 2,$ these cases are impossible. Therefore, assume that 
$k\geq 3.$ Then, $n\geq 5.$ Let $\alpha $ be positive real root of $\Psi
_{k}(x)$ given in (\ref{0.2}). Then $2<\alpha <\varphi ^{2}<3$ by Lemma \ref%
{L1} (b). Besides, it is seen that $10^{m-1}<P_{n}^{(k)}<10^{m}.$ Thus,
using the inequality (\ref{1.3}), we get 
\begin{equation*}
(n-2)\frac{\log 2}{\log 10}<m<(n-1)\frac{\log 3}{\log 10}+1,
\end{equation*}%
which implies that 
\begin{equation}
\frac{3n}{20}<m<\frac{3n}{4}  \label{3.1}
\end{equation}%
for $n\geq 5.$ Now, rearranging the equation (\ref{1.1}) as%
\begin{equation*}
P_{n}^{(k)}-g_{k}(\alpha )\alpha ^{n}+\dfrac{d}{9}=d\dfrac{10^{m}}{9}%
-g_{k}(\alpha )\alpha ^{n}
\end{equation*}%
and taking absolute value of both sides, we get 
\begin{equation}
\left\vert d\dfrac{10^{m}}{9}-g_{k}(\alpha )\alpha ^{n}\right\vert <\frac{3}{%
2}  \label{3.2}
\end{equation}%
using the inequality (\ref{w}). If we divide both sides of the inequality (%
\ref{3.2}) by $g_{k}(\alpha )\alpha ^{n}$, from Lemma \ref{L1}, we get 
\begin{equation}
\left\vert 10^{m}\alpha ^{-n}\frac{d(g_{k}(\alpha ))^{-1}}{9}-1\right\vert <%
\frac{3}{2g_{k}(\alpha )\alpha ^{n}}<\frac{1}{0.552\cdot \alpha ^{n}}<\frac{%
5.5}{\alpha ^{n}}.  \label{3.3}
\end{equation}%
In order to use the result of Matveev Theorem \ref{T2}, we take 
\begin{equation*}
\left( \gamma _{1},b_{1}\right) :=\left( 10,m\right) ,~\left( \gamma
_{2},b_{2}\right) :=\left( \alpha ,-n\right) ,~\left( \gamma
_{3},b_{3}\right) :=\left( \frac{9\cdot g_{k}(\alpha )}{d},-1\right) .
\end{equation*}%
The number field containing $\gamma _{1},~\gamma _{2},$ and $\gamma _{3}$
are $\mathbb{K}=\mathbb{Q}(\sqrt{\alpha }),$ which has degree $D=k.$ We show
that the number%
\begin{equation*}
\Lambda _{1}:=10^{m}\alpha ^{-n}\frac{d(g_{k}(\alpha ))^{-1}}{9}-1
\end{equation*}%
is nonzero. Contrast to this, assume that $\Lambda _{1}=0$. Then 
\begin{equation*}
d\dfrac{10^{m}}{9}=\alpha ^{n}g_{k}(\alpha )=\frac{\alpha -1}{(k+1)\alpha
^{2}-3k\alpha +k-1}\alpha ^{n}.
\end{equation*}%
Conjugating the above equality by some automorphisim of the Galois group of
the splitting field of $\Psi _{k}(x)$ over $%
\mathbb{Q}
$ and taking absolute values, we get 
\begin{equation*}
d\dfrac{10^{m}}{9}=\left\vert \frac{\alpha _{i}-1}{(k+1)\alpha
_{i}^{2}-3k\alpha _{i}+k-1}\alpha _{i}^{n}\right\vert
\end{equation*}%
for some $i>1,$ where $\alpha =\alpha _{1},\alpha _{2},\ldots ,\alpha _{k}$
are the roots of $\Psi _{k}(x).$ Using (\ref{1.4}) and that $|\alpha _{i}|<1$%
, we obtain from the last equality that 
\begin{eqnarray*}
d\dfrac{10^{m}}{9} &=&\left\vert \frac{\alpha _{i}^{k}-\alpha _{i}^{k-1}}{%
\alpha _{i}^{k+1}-\alpha _{i}^{k-1}-k}\right\vert \left\vert \alpha
_{i}\right\vert ^{n} \\
&<&2,
\end{eqnarray*}%
which is impossible since $m\geq 2.$ Therefore $\Lambda _{1}\neq 0.$
Moreover, since $h(10)=\log 10,$ $h(\gamma _{2})=\dfrac{\log \alpha }{k}<%
\dfrac{\log 3}{k}$ by (\ref{2.1}) and%
\begin{eqnarray*}
h(\gamma _{3}) &=&h(\frac{9\cdot g_{k}(\alpha )}{d})\leq
h(9)+h(9)+h(g_{k}(\alpha )) \\
&\leq &\log 81+4\log k\leq 8\log k
\end{eqnarray*}%
by (\ref{2.5}), we can take $A_{1}:=k\log 2,~A_{2}:=\log 3,$ and $%
A_{3}:=8k\log k.$ Also, since $m\leq 3n/4,$ it follows that $B:=n.$ Thus,
taking into account the inequality (\ref{3.3}) and using Theorem \ref{T2},
we obtain{\small 
\begin{equation*}
\frac{5.5}{\alpha ^{n}}>\left\vert \Lambda _{1}\right\vert >\exp \left(
-1.4\cdot 30^{6}\cdot 3^{4.5}\cdot k^{2}(1+\log k)(1+\log n)\left( k\log
10\right) \left( \log 3\right) \left( 8k\log k\right) \right)
\end{equation*}%
}and so 
\begin{equation*}
n\log \alpha -\log (5.5)<1.4\cdot 30^{6}\cdot 3^{4.5}\cdot k^{2}(1+\log
k)(1+\log n)\left( k\log 10\right) \left( \log 3\right) \left( 8k\log
k\right) ,
\end{equation*}%
where we have used the fact that $1+\log y<2\log y$ for all $y\geq 3.$ From
the last inequality, a quick computation with Mathematica yields to\textit{\ 
}%
\begin{equation*}
n\log \alpha <1.16\cdot 10^{13}\cdot k^{4}\cdot (\log k)^{2}\cdot \log n
\end{equation*}%
or \textit{\ }%
\begin{equation}
n<1.68\cdot 10^{13}\cdot k^{4}\cdot (\log k)^{2}\cdot \log n.  \label{3.5}
\end{equation}%
The inequalitry (\ref{3.5}) can be rearranged as 
\begin{equation*}
\frac{n}{\log n}<1.68\cdot 10^{13}\cdot k^{4}\cdot (\log k)^{2}.
\end{equation*}%
Using the fact that 
\begin{equation*}
\text{if }A\geq 3\text{ and }\frac{n}{\log n}<A,\text{ then }n<2A\log A,
\end{equation*}%
we obtain 
\begin{eqnarray}
n &<&3.36\cdot 10^{13}\cdot k^{4}\cdot (\log k)^{2}\log \left( 1.68\cdot
10^{13}\cdot k^{4}\cdot (\log k)^{2}\right)  \label{3.6} \\
&<&3.36\cdot 10^{12}\cdot k^{4}\cdot (\log k)^{2}(30.5+4\log k+2\log (\log
k))  \notag \\
&<&3.36\cdot 10^{12}\cdot k^{4}\cdot (\log k)^{2}(34\log k)  \notag \\
&<&1.15\cdot 10^{15}\cdot k^{4}\cdot (\log k)^{3},  \notag
\end{eqnarray}%
where we have used the fact that $30.5+4\log k+2\log (\log k)<34\log k$ for
all $k\geq 3.$

Let $k\in \lbrack 3,400].$ Now, let us try to reduce the upper bound on $n$
applying Lemma \ref{L2}. Let 
\begin{equation*}
z_{1}:=m\log 10-n\log \alpha +\log \left[ \frac{d}{9}\left( g_{k}(\alpha
)\right) ^{-1}\right] .
\end{equation*}

\bigskip and $x:=e^{z_{1}}-1$. Then from (\ref{3.3}), we get 
\begin{equation*}
\left\vert x\right\vert =\left\vert e^{z_{1}}-1\right\vert <\frac{5.5}{%
\alpha ^{n}}<0.2
\end{equation*}%
for $n\geq 5.$ Choosing $a:=$ $0.2,$ we obtain the inequality

\begin{equation*}
|z_{1}|=\left\vert \log (x+1)\right\vert <\frac{\log (10/8)}{(0.2)}\cdot 
\frac{5.5}{\alpha ^{n}}<\frac{6.14}{\alpha ^{n}}
\end{equation*}%
by Lemma \ref{L3}. Thus, it follows that 
\begin{equation*}
0<\left\vert m\log 10-n\log \alpha +\log \left[ \frac{d}{9}\left(
g_{k}(\alpha )\right) ^{-1}\right] \right\vert <\frac{6.14}{\alpha ^{n}}.
\end{equation*}%
Dividing this inequality by $\log \alpha ,$ we get 
\begin{equation}
0<|m\gamma -n+\mu |<A\cdot B^{-w},  \label{3.7}
\end{equation}%
where 
\begin{equation*}
\gamma :=\dfrac{\log 10}{\log \alpha }\notin 
\mathbb{Q}
,~\mu :=\dfrac{\log \left( \frac{d}{9}\left( g_{k}(\alpha )\right)
^{-1}\right) }{\log \alpha },~A:=8.86,~B:=\alpha \text{, and }w:=n.
\end{equation*}%
If we take 
\begin{equation*}
M:=\left\lfloor 1.15\cdot 10^{15}\cdot k^{4}\cdot (\log k)^{3}\right\rfloor ,
\end{equation*}%
which is an upper bound on $m$ since $m<n<1.15\cdot 10^{15}\cdot k^{4}\cdot
(\log k)^{3}$ by (\ref{3.6}), we found that $q_{71},$ the denominator of the 
$71$st convergent of $\gamma $ exceeds $6M.$ Furthermore, a quick
computation with Mathematica gives us that the value%
\begin{equation*}
\dfrac{\log \left( Aq_{71}/\epsilon \right) }{\log B}
\end{equation*}%
is less than $99.3.$ So, if the inequality (\ref{3.7}) has a solution, then 
\textit{\ } 
\begin{equation*}
n<\dfrac{\log \left( Aq_{71}/\epsilon \right) }{\log B}\leq 99.3,
\end{equation*}%
that is, $n\leq 99.$ In this case, $m<75$ by (\ref{3.1}). A quick
computation with Mathematica gives us that the equation $P_{n}^{(k)}=d\left( 
\frac{10^{m}-1}{9}\right) $ has the solutions for $%
(n,k,d,m)=(5,3,3,2),(6,4,8,2)$ in the intervals $n\in \left[ 5,99\right]
,~m\in (2,75)$ and $k\in \left[ 3,400\right] .$ Thus, this completes the
analysis in the case $k\in \left[ 3,400\right] .$

From now on, we can assume that $k>400.$ Then we can see from (\ref{3.6})
that the inequality%
\begin{equation}
n<1.15\cdot 10^{15}\cdot k^{4}\cdot (\log k)^{3}<\varphi ^{k/2}  \label{x}
\end{equation}%
holds for $k>400.$

By Lemma $7$ given in \cite{zfr2}, we have 
\begin{equation}
g_{k}(\alpha )\alpha ^{n}=\frac{\varphi ^{2n}}{\varphi +2}+\frac{\delta }{%
\varphi +2}+\eta \varphi ^{2n}+\eta \delta ,  \label{3.14}
\end{equation}%
where 
\begin{equation}
\left\vert \delta \right\vert <\dfrac{\varphi ^{2n}}{\varphi ^{k/2}}\text{
and }\left\vert \eta \right\vert <\frac{3k/2}{\varphi ^{k}}.  \label{3.13}
\end{equation}%
So, using (\ref{3.2}), (\ref{3.14}) and (\ref{3.13}), we obtain 
\begin{eqnarray}
\left\vert d\dfrac{10^{m}}{9}-\frac{\varphi ^{2n}}{\varphi +2}\right\vert
&=&\left\vert \left( d\dfrac{10^{m}}{9}-g_{k}(\alpha )\alpha ^{n}\right) +%
\frac{\delta }{\varphi +2}+\eta \varphi ^{2n}+\eta \delta \right\vert
\label{3.15} \\
&\leq &\left\vert d\dfrac{10^{m}}{9}-g_{k}(\alpha )\alpha ^{n}\right\vert +%
\frac{\left\vert \delta \right\vert }{\varphi +2}+\left\vert \eta
\right\vert \varphi ^{2n}+\left\vert \eta \right\vert \left\vert \delta
\right\vert  \notag \\
&<&\frac{3}{2}+\dfrac{\varphi ^{2n}}{\varphi ^{k/2}\left( \varphi +2\right) }%
+\frac{3k\varphi ^{2n}}{2\varphi ^{k}}+\dfrac{3k\varphi ^{2n}}{2\varphi
^{3k/2}}.  \notag
\end{eqnarray}%
Dividing both sides of the above inequality by $\dfrac{\varphi ^{2n}}{%
\varphi +2},$ we get{\small 
\begin{eqnarray}
\left\vert 10^{m}\varphi ^{-2n}\frac{d}{9}\left( \varphi +2\right)
-1\right\vert &<&\frac{3\left( \varphi +2\right) }{2\varphi ^{2n}}+\dfrac{1}{%
\varphi ^{k/2}}+\frac{3k\left( \varphi +2\right) }{2\varphi ^{k}}+\dfrac{%
3k\left( \varphi +2\right) }{2\varphi ^{3k/2}}  \label{3.16} \\
&<&\dfrac{0.15}{\varphi ^{k/2}}+\dfrac{1}{\varphi ^{k/2}}+\dfrac{0.005}{%
\varphi ^{k/2}}+\dfrac{0.005}{\varphi ^{k/2}}=\dfrac{1.16}{\varphi ^{k/2}}, 
\notag
\end{eqnarray}%
} where we have used the facts that 
\begin{equation*}
\frac{3k\left( \varphi +2\right) }{2\varphi ^{k}}<\dfrac{0.005}{\varphi
^{k/2}}\text{ and }\dfrac{3k\left( \varphi +2\right) }{2\varphi ^{3k/2}}<%
\dfrac{0.005}{\varphi ^{k/2}}\text{ for }k>400.
\end{equation*}%
In order to use the result of Matveev Theorem \ref{T2}, we take 
\begin{equation*}
\left( \gamma _{1},b_{1}\right) :=\left( 10,m\right) ,~\left( \gamma
_{2},b_{2}\right) :=\left( \varphi ,-2n\right) ,~\left( \gamma
_{3},b_{3}\right) :=\left( \frac{d\left( \varphi +2\right) }{9},1\right) .
\end{equation*}%
The number field containing $\gamma _{1},~\gamma _{2},$ and $\gamma _{3}$
are $\mathbb{K}=\mathbb{Q}(\sqrt{5}),$ which has degree $D=2.$ We show that
the number 
\begin{equation*}
\Lambda _{1}:=10^{m}\varphi ^{-2n}\frac{d}{9}\left( \varphi +2\right) -1
\end{equation*}
is nonzero. Contrast to this, assume that $\Lambda _{1}=0.$ Then $10^{m}%
\frac{d}{9}\left( \varphi +2\right) =\varphi ^{2n}$ and conjugating this
relation in $\mathbb{Q}(\sqrt{5}),$ we get $10^{m}\frac{d}{9}\left( \beta
+2\right) =\beta ^{2n},$ where $\beta =\frac{1-\sqrt{5}}{2}=\overline{%
\varphi }.$ So, we have%
\begin{equation*}
\dfrac{\varphi ^{2n}}{\varphi +2}=\dfrac{\beta ^{2n}}{\beta +2},
\end{equation*}%
which implies that 
\begin{equation*}
\dfrac{\varphi ^{4n}}{\varphi +2}=\dfrac{1}{\beta +2}<1.
\end{equation*}
The last inequality is impossible for $n\geq 5.$ Therefore $\Lambda _{1}\neq
0.$ Moreover, since 
\begin{equation*}
h(\gamma _{1})=h(10)=\log 10,h(\gamma _{2})=h(\varphi )\leq \frac{\log
\varphi }{2}
\end{equation*}%
and 
\begin{equation*}
h(~\gamma _{3})\leq h(9)+h(d)+h(\varphi )+h(2)+\log 2\leq \log 324+\frac{%
\log \varphi }{2}
\end{equation*}%
by (\ref{2.3}), we can take $A_{1}:=2\log 10,~A_{2}:=\log \varphi $, and $%
A_{3}:=\log \left( 324^{2}\varphi \right) .$ Also, since $m<3n/4,$ we can
take $B:=2n.$ Thus, taking into account the inequality (\ref{3.16}) and
using Theorem \ref{T2}, we obtain%
\begin{equation*}
\left( 1.16\right) \cdot \varphi ^{-k/2}>\left\vert \Lambda _{1}\right\vert
>\exp \left( C\cdot (1+\log 2n)\left( 2\log 10\right) \left( \log \varphi
\right) \log \left( 324^{2}\varphi \right) \right) ,
\end{equation*}%
where $C=-1.4\cdot 30^{6}\cdot 3^{4.5}\cdot 2^{2}\cdot (1+\log 2)$. This
implies that 
\begin{equation*}
\frac{k}{2}\log \varphi -\log (1.16)<2.59\cdot 10^{13}\cdot (1+\log 2n)
\end{equation*}%
or 
\begin{equation}
k<2.16\cdot 10^{14}\cdot \log 2n,  \label{3.17}
\end{equation}%
where we have used the fact that $(1+\log 2n)<2\log 2n$ for $n\geq 5.$ On
the other hand, from (\ref{x}), we get 
\begin{eqnarray*}
\log 2n &<&\log \left( 2.3\cdot 10^{15}\cdot k^{4}\cdot (\log k)^{3}\right)
\\
&<&35.4+4\log k+3\log (\log k) \\
&<&40\log k.
\end{eqnarray*}%
So, from (\ref{3.17}), we obtain 
\begin{equation*}
k<2.16\cdot 10^{14}\cdot 40\log k,
\end{equation*}%
which implies that 
\begin{equation}
k<3.5\cdot 10^{17}.  \label{3.18}
\end{equation}%
To reduce this bound on $k$, we use Lemma \ref{L2}. Substituting this bound
of $k$ into (\ref{x}), we get $n<1.14\cdot 10^{90},$ which implies that $%
m<8.55\cdot 10^{89}$ by (\ref{3.1}).

Now, let 
\begin{equation*}
z_{2}:=m\log 10-2n\log \varphi +\log \left( \frac{d}{9}\left( \varphi
+2\right) \right) .
\end{equation*}%
and $x:=1-e^{z_{2}}$. Then 
\begin{equation*}
|x|=\left\vert 1-e^{z_{2}}\right\vert <\dfrac{1.16}{\varphi ^{k/2}}<0.6
\end{equation*}%
by (\ref{3.16}). Choosing $a:=$ $0.6,$ obtaint the inequality

\begin{equation*}
|z_{2}|=\left\vert \log (x+1)\right\vert <\frac{\log (5/2)}{0.6}\cdot \dfrac{%
1.16}{\varphi ^{k/2}}<\frac{1.78}{\varphi ^{k/2}}
\end{equation*}%
by Lemma \ref{L3}. That is, 
\begin{equation*}
0<\left\vert m\log 10-2n\log \varphi +\log \left( \frac{d}{9}\left( \varphi
+2\right) \right) \right\vert <\frac{1.78}{\varphi ^{k/2}}.
\end{equation*}%
Dividing both sides of the above inequality by $\log \varphi ,$ we get 
\begin{equation}
0<|m\gamma -2n+\mu |<A\cdot B^{-w},  \label{3.20}
\end{equation}%
where 
\begin{equation*}
\gamma :=\dfrac{\log 10}{\log \varphi }\notin 
\mathbb{Q}
,~\mu :=\dfrac{\log \left( \frac{d}{9}\left( \varphi +2\right) \right) }{%
\log \varphi },~A:=3.7,~B:=\varphi \text{, and }w:=k/2.
\end{equation*}%
If we take $M:=8.55\cdot 10^{89}$, which is an upper bound on $m$, we found
that $q_{180},$ the denominator of the $180$ th convergent of $\gamma $
exceeds $6M.$ Furthermore, a quick computation with Mathematica gives us
tthat he value%
\begin{equation*}
\dfrac{\log \left( Aq_{180}/\epsilon \right) }{\log B}
\end{equation*}%
is less than $444.7.$ So, if the inequality (\ref{3.20}) has a solution,
then \textit{\ } 
\begin{equation*}
\frac{k}{2}<\dfrac{\log \left( Aq_{180}/\epsilon \right) }{\log B}\leq 444.7,
\end{equation*}%
that is, $k\leq 889.$ Hence, from (\ref{x}), we get $n<2.25\cdot 10^{29},$
which implies that $m<1.69\cdot 10^{29}$ since $m<3n/4$ by (\ref{3.1}). If
we apply Lemma \ref{L2} again with $M:=1.69\cdot 10^{29},$ we found that $%
q_{60},$ the denominator of the $60$ th convergent of $\gamma $ exceeds $6M.$
After doing this, then a quick computation with Mathematica show that the
inequality (\ref{3.20}) has solution only for $k<313.$ This contradicts the
fact that $k>400.$ This completes the proof.%
\endproof%

\end{document}